\def\timestamp{%
Time-stamp: <embbool.tex: Monday 03-12-2001 at 15:25:42 (cet)>}
\def\stripname Time-stamp: <#1 #2>{#2}
\edef\filedate{\expandafter\stripname\timestamp}

\documentclass[a4paper]{amsart}[2000/06/02]
\usepackage{eucal,amsfonts}
\newcommand\calA{\mathcal{A}}
\newcommand\calB{\mathcal{B}}
\newcommand\calC{\mathcal{C}}

\newcommand\calF{\mathcal{F}}
\newcommand\calI{\mathcal{I}}

\newcommand\calP{\mathcal{P}}
\newcommand\calR{\mathcal{R}}

\newcommand\calX{\mathcal{X}}
\newcommand\calY{\mathcal{Y}}
\renewcommand\newsymbol[5]{%
\DeclareMathSymbol#1{#3}{\ifcase #2\or AMSa\or AMSb\fi}{"#4#5}}
\newcommand\newbbbletter[2]{%
\DeclareMathSymbol#1{0}{AMSb}{`#2}}
\let\emptyset \undefined
\let\ge       \undefined
\let\le       \undefined
\def\c{\mathfrak c}
\def\A{\mathbb A}
\def\B{\mathbb B}
\def\C{\mathbb C}
\def\P{\mathbb P}
\def\Q{\mathbb Q}
\def\M{\mathbb M}
\newcommand\fin{\operatorname{fin}}
\def\PNfin{\calP( \N) /\fin}
\def\MASC{\mathrm{MA}_{\sigma\mathrm{-centered}}}
\def\PFA{\mathrm{PFA}}
\def\MASL{\mathrm{MA}_{\sigma\mathrm{-linked}}}
\def\<{\langle}
\def\>{\rangle}
\def\CH{\mathrm{CH}}

\newsymbol\restr       1216
\newsymbol\mapdiagbin  124D    
\newsymbol\le          1336    \let\leq\le
\newsymbol\ge          133E    \let\geq\ge
\newsymbol\emptyset    203F
\newsymbol\notle       230A
\newsymbol\subsetneq   2328
\newsymbol\nsubseteq   232A
\newbbbletter\Cantor   C
\newbbbletter\I        I
\let\L\undefined
\newbbbletter\L        L
\newbbbletter\N        N
\newbbbletter\psarc    P
\newbbbletter\R        R
\newcommand\cl[1]{\operatorname{cl}#1}

\newcommand\St{\operatorname{St}}
\newcommand\dom{\operatorname{dom}}
\newcommand\rng{\operatorname{rng}}
\newcommand\Aut{\operatorname{Aut}}
\newcommand\CO{\operatorname{CO}}
\newcommand\RO{\operatorname{RO}}

\newcommand\id{\operatorname{id}}

\DeclareMathSymbol\mapdiag{1}{symbols}{"34}
\let\epsilon\varepsilon
\let\meet\wedge

\newcommand\0{\mathbf{0}}
\newcommand\1{\mathbf{1}}
\DeclareMathSymbol{\lor}{\mathrel}{symbols}{"5F}
\DeclareMathSymbol{\land}{\mathrel}{symbols}{"5E}
\def\nfrac#1/#2{\leavevmode\kern.1em
  \raise.5ex\hbox{\the\scriptfont0 #1}\kern-.1em
  /\kern-.15em\lower.25ex\hbox{\the\scriptfont0 #2}}

\newtheorem{theorem}{Theorem}[section]
\newtheorem{corollary}[theorem]{Corollary}
\newtheorem{lemma}[theorem]{Lemma}
\newtheorem{question}[theorem]{Question}
\newtheorem{proposition}[theorem]{Proposition}

\theoremstyle{definition}
\newtheorem{defn}[theorem]{Definition}
\numberwithin{equation}{section}
\theoremstyle{remark}

\usepackage{epsfig}
\begin{document}

\title[Embeddings into $\PNfin$ and extension of automorphisms]%
      {Embeddings into $\PNfin$ and extension of automorphisms}

\author{A. Bella}
\address[A. Bella, P. Ursino]{Dipartimento di Matematica\\
              Universit\`a di Catania\\
              95125 Catania\\
              Italy}
\email{bella@dmi.unict.it, ursino@dmi.unict.it}

\author{A. Dow}
\address[A. Dow]{Mathematics Department\\ 
         UNC Charlotte\\
         9201 University City Blvd.\\ 
         Charlotte, NC 28223\\
         USA}
\email{adow@uncc.edu}

\author{K. P. Hart}
\address[K. P. Hart]{Faculty of Information Technology and Systems\\TU Delft\\
         Postbus 5031\\2600~GA {} Delft\\the Netherlands}
\email{k.p.hart@its.tudelft.nl}
%\urladdr{http://aw.twi.tudelft.nl/\~{}hart}

\author[M. Hru\v{s}\'ak]{M. Hru\v{s}\'ak${}^\dag$}
\thanks{\leavevmode\llap{${}^\dag$}%
        The research of the fourth author was supported in part by
        the Netherlands Organization for Scientific Research (NWO) ---
        Grant 613.007.039, and in part by the Grant Agency of the Czech Republic --- Grant GA\v CR 201/00/1466.}
\author{J. van Mill}
\address[M. Hru\v{s}\'ak, J. van Mill]{Faculty of Sciences\\
         Division of Mathematics\\
         Vrije Universiteit\\
         De Boelelaan 1081\textsuperscript{a}\\
         1081 HV Amsterdam\\
         The Netherlands}
\email{michael@cs.vu.nl, j.van.mill@cs.vu.nl}
\author{P. Ursino}
\date{\filedate}

\begin{abstract} Given a Boolean algebra $\B$ and an embedding $e:\B\to \PNfin$ 
we consider the possibility of extending each or some
automorphism of $\B$ to the whole $\PNfin$. Among other
things, we show, assuming~$\CH$, that for a wide class of Boolean algebras there are embeddings for which no non-trivial
automorphism can be extended.
\end{abstract}



\maketitle

\section{Introduction}

A general problem in the theory of Boolean algebras concerns the
possibility of extending all or some automorphism  with respect to a
fixed embedding \cite {St}. More precisely, given two Boolean
algebras $\A$ and $\B$  one asks whether there is an
embedding $e: \B \to \A$ for which   each automorphism of
$\B$ can be extended to  an automorphism of $\A$. In
the opposite direction, one may wonder if there is an embedding
for which no non-trivial automorphism of $\B$ extends.

For instance, P. \v St\v ep\' anek in \cite{St} quotes as open the problem
whether any Boolean algebra $\B$ can be embedded into a
homogeneous Boolean algebra $\A$ in such a way that no non-trivial 
automorphism of $\B$ extends to an automorphism of $\A$.
A result in this direction due to S. Koppelberg says that, assuming
$\diamondsuit$, there exists an $\aleph_1$-Suslin tree $T$ such
that the corresponding complete Boolean algebra $\A(T)$ is
homogeneous and  has  a regular complete subalgebra  onto which
no non-trivial automorphism of $\A(T)$ restricts.

This paper focuses  on the case  $\A=\PNfin$. The main reason for this is
Parovi\v cenko's theorem:

\begin{theorem} Every Boolean algebra of size $\leq\omega_1$ embeds into $\PNfin$.
\end{theorem}

So assuming the Continuum Hypothesis, for algebras of size at most continuum
at least the minimal requirement of an existence of an embedding is satisfied.
One should also mention the following result of van Douwen and Przymusi\' nski 
(see \cite{vDP} or \cite{vM}) 

\begin{theorem}[MA] Every Boolean algebra of size $<\c$ embeds into $\PNfin$.
\end{theorem}

However, this result 
does not hold for all algebras of size $\c$. In fact,
assuming the Proper Forcing Axiom even algebras as nice as the measure algebra $\M=\text{Borel}(\R)/\text{Null}$
do not embed into $\PNfin$ (see \cite{DH}). \par
The investigation of how the measure algebra can be embedded into $\PNfin$ was suggested 
by the last mentioned author as an 
attempt to compare, as composition groups, certain measure-preserving  endomorphisms 
on the real line \cite{GPU}  with the permutations of the integers.
Since a straightforward attempt presents unsurmountable obstacles, 
a possible way to approach this problem indirectly could be to view both
of them as embedded in the larger group $\Aut(\PNfin)$.
Indeed, the latter give rise to automorphisms of $\PNfin$, 
and as we show, the same happens for the former.\par

The results presented here show that, under $\CH$, for a wide class of Boolean algebras, including the measure algebra
there are embeddings such that all automorphisms
can be extended, as well as  embeddings such that no non-trivial automorphism can be
extended.\par
All basic concepts of Boolean algebra and notation  can be found in
\cite{Ko}, set theoretic notation follows \cite{Ku}. Henceforth, $\CH$ stands for the Continuum Hypothesis, 
MA, $\MASC$ and $\MASL$ indicate
the use of Martin's Axiom or its variants, while $\PFA$ stands for the Proper Forcing Axiom.
As usual, given a Tychonoff space $X$, $\beta X$
denotes the \v Cech-Stone compactification of $X$, and  $X^*$
the  \v Cech-Stone remainder $\beta X\setminus X$. $\CO(X)$ and $\RO(X)$
are the Boolean algebras of clopen and regular open subsets of~$X$, 
respectively. 
In particular, $\CO(2^\kappa)$ denotes free Boolean algebra on $\kappa$~many 
generators. The Stone space of a Boolean algebra~$\B$ is written $\St(\B)$.
\section{Nice embeddings}

The following fact that, assuming $\CH$, nice embeddings do exist appeared in 
\cite{Gr} with a different proof. 

\begin{theorem}[$\CH$]\label{grzech}
 For every Boolean algebra $\mathbb B$ of size at most continuum there is an embedding $e:\mathbb B\to \PNfin$
such that every automorphism of $\mathbb B$ can be extended to an automorphism of $\PNfin$.

\end{theorem}

\begin{proof} Given a Boolean algebra $\mathbb B$ of size at most $\c$ 
note that the algebra $\mathbb B^\omega/\text{fin}$ is isomorphic to $\PNfin$
(\cite{vM}).
Define an embedding $e:\B\to\B^\omega/\text{fin}$ by $e(B)=[ (B,B,B,\dots)]$. Clearly, if $h$ is an automorphism of 
$\B$, putting
$H([(A_n)_{n\in\omega}]) =[h(A_n)_{n\in\omega}]$ defines an extension of~$h$,
which is an automorphism of $\B^\omega/\text{fin}$.
\end{proof}

The above theorem can be strengthened by requiring that the embedding $e$ be such that every automorphism
$h$ has the largest possible number, i.e., $2^\c$ many, extensions. This is accomplished by identifying 
$\PNfin$ with $\A^\omega/\text{fin}$, where $\A = \B \oplus \CO(2^\c)$.\par
A natural question as to whether the hypothesis of Theorem \ref{grzech} can be weakened was also addressed in \cite{Gr}
where it is shown that

\begin{theorem} It is relatively consistent with $\MASL$ that for every Boolean algebra 
$\mathbb B$ of size at most continuum there is an embedding $e:\mathbb B\to \PNfin$
such that each automorphism of $\mathbb B$ can be extended to an automorphism of $\PNfin$.
\end{theorem}

It should be noted here that for certain Boolean algebras $\B$ such as any countable algebra all embeddings are good.
Let $e:\B\to\PNfin$ be an embedding. We say that \emph{$e$ lifts} if there is a Boolean algebra homomorphism 
$E:\B\to\calP(\N)$ such that $E(B)\in e(B)$ for every $B\in\B$.

\begin{theorem}[$\MASC(\kappa)$] \label{goodemb_MA}
Let $\B$ be a Boolean algebra of size $\kappa$ and let $e:\B\to\PNfin$ be an embedding which lifts. 
Then every automorphism of $\B$ extends to an automorphism of $\PNfin $.
\end{theorem}

\begin{proof} Let a Boolean algebra $\B$ and an embedding $e:\B\to \PNfin$ be given.
Let $E:\B\to\calP(\N)$ be a lifting of $e$. Let $h$ be an automorphism of $\B$.
 Consider the following partial order
$$\P =\{ \<s,F\>: s\text{ is a finite partial one-to-one function }\N\to\N,\  F\in [\B]^{<\aleph_0}\}$$
ordered by $\<s,F\>\leq\<t,G\>$ if $s\supseteq t$, $F\supseteq G$ and for every $n\in\dom(s)\setminus\dom(t)$ and every $B\in G$,
$n\in E(B)$ if and only if $s(n)\in E(h(B))$.\par

The partial order $\P$ is $\sigma$-centered as conditions with the same working part are compatible. It is also easily
seen that the following sets are dense in $\P$:
\begin{itemize}
\item $D_n=\{ \<s,F\>\in\P: n\in\dom(s)\}$,
\item $R_n=\{ \<s,F\>\in\P: n\in\rng(s)\}$,
\item $E_B=\{ \<s,F\>\in\P: B\in F\}$.
\end{itemize}

By $\MASC$ there is a filter $G$ on $\P$ which intersects all of these dense sets. Let 
$\pi=\bigcup\{s:(\exists F\in [\B]^{<\aleph_0})(\<s,F\>\in G)\}$. 
It is straightforward to check that $\pi$~is a permutation on~$\N$ 
that defines an automorphism of $\PNfin$ extending~$h$.
\end{proof}  

\begin{corollary}[$\MASC(\kappa)$] Every embedding $e:\CO(2^{\kappa})\to\PNfin$
is such that every automorphism of $\CO(2^{\kappa})$ extends to an automorphism of $\PNfin$.
\end{corollary}

\begin{proof} It is easy to see that every $e:\CO(2^{\kappa})\to\PNfin$ lifts. 
\end{proof}

Note the fundamental difference between the two proofs presented so far; 
whereas in Theorem \ref{grzech} it was in effect shown that the map that sends
 $h$ to~$H$   is actually an injective group homomorphism from $\Aut(\B)$ 
into~$\Aut(\PNfin)$, the method of the proof of 
Theorem \ref{goodemb_MA} does not seems to produce such an embedding.

\section{Ugly Embeddings}

This section is devoted to showing that, assuming $\CH$, many Boolean algebras can be embedded into $\PNfin$ in such a way that only the
trivial automorphism extends to an automorphism of $\PNfin$. 
Recall that if $A$ and~$B$ are elements of a Boolean algebra $\B$ then 
\emph{$A$ splits $B$} if both $A\meet B$ and $A^c\meet B$ are non-zero.
Denote by $\mathfrak u(\B)$ the minimal character of an ultrafilter on $\B$.

\begin{lemma}[$\CH$] \label{split_emb}
Let $\B$ be a Boolean algebra of size $\c$ such that $\mathfrak u(\B)>\omega$. 
Then for every family of ultrafilters $\{p_\alpha  :\alpha  <\c\}\subseteq \St(\B)$ 
 there exists an embedding $e: \B \to  \PNfin$ and a  set
$\{C_\alpha  : \alpha <\c\}\subseteq \PNfin^+$ 
such that
\begin{enumerate}
\item $C_\alpha <e(B)$ for each $B\in p_\alpha $ and $\alpha <\c$,
\item If $A\in\PNfin^+$ and $A\meet C_\alpha =\0$ for every $\alpha <\c$ then 
there is a  $B\in \B$ such that $e(B)$ splits $A$.
\end{enumerate}
\end{lemma}

\begin{proof} 
Enumerate  $\B$ as $\{B_{\alpha} : \alpha\in\omega_1\}$ and $\PNfin$ as $\{A_{\alpha} : \alpha\in\omega_1\}$. Recursively construct
an increasing chain  $\{\B_\beta :\beta\in \omega_1 \}$ of countable subalgebras of $\B$ and compatible
embeddings $e_\beta :\B_\beta \to \PNfin$ together with a set of non-zero pairwise disjoint elements
$\{C_\beta : \beta <\omega_1\}\subseteq \PNfin$ so that for $\alpha<\c$:
\begin{enumerate}
\item $B_\alpha\in \B_\alpha$,
\item $C_\gamma  < e_\alpha(B)$ for every  $B\in p_\gamma \cap \B_\alpha$ and $\gamma\le \alpha$,
\item if $(\forall\gamma \le \alpha)(A_\alpha \meet C_\gamma =\0) $  then  $(\exists B\in
\B_\alpha)(e_\alpha(B)$ splits $A_\alpha)$.
\end{enumerate}
At stage $\alpha<\c$ we have to   define the subalgebra $\B_\alpha$, the embedding $e_\alpha$ and the element $C_\alpha\in\PNfin$.
Let $\B' =\bigcup \{\B_\beta :\beta <\alpha \}$ and let $e'=\bigcup\{e_\beta : \beta <\alpha\}$ be the corresponding
embedding. Let $\C $ be   the subalgebra of $\B$ generated by $\B'\cup \{B_\alpha \}$. Let 
 $c:\C \to \PNfin$ be an embedding which extends $e'$ 
and such that $C_\beta < c (B_\alpha )$ whenever $B_\alpha\in p_\beta$, and $C_\beta \meet c (B_\alpha )=\0$ if $B_\alpha\not\in p_\beta$. 
Then pick a $C_\alpha \in \PNfin^+$ so that  $C_\alpha <c(B)$  for all $B\in p_\alpha \cap
\C $ and $C_\alpha \meet C_\beta =\0$ for each $\beta<\alpha $.\par      

If there is an element $B\in\C $ such that $c(B)$ splits  $A_\alpha $ or if $A_\alpha
\meet C_\beta >\0$ for some $\beta \le \alpha $, then set 
$\B_\alpha =\C $ and $e_\alpha =c$. If not, let $\calC=\{ B\in  
\B'_{\alpha}: A_\alpha <c(B)\}$. The set  $\calC$ generates a
filter in $\B$. Being countably generated, this filter cannot
be an ultrafilter so there is a $G\in\B$ which splits every element of $\calC$.
Let $\B_\alpha $ be the subalgebra of $\B$
generated by $\C \cup\{G\}$. All that has to be done is to extend $c$ to an embedding $e_\alpha :\B_\alpha  \to \PNfin$ 
so that  $e_\alpha (G)$ splits $A_\alpha $. To that end let
$$F=\{c(B):B<G\}\cup\{ C_\beta:\beta\leq\alpha \ \& \ G\in x_\beta\},$$
$$G=\{c(B):G<B\}\cup\{ {C_\beta}^c:\beta\leq\alpha \ \& \ G\not\in x_\beta\},$$
$$H=\{c(B):B\not<G \ \& \ B\not>G\}\cup\{A_\alpha\}.$$
As $\PNfin$ satisfies condition $\calR_\omega$ (see \cite{vM}) there is an element $Y\in\PNfin$ which separates
$F$ from $G$ and splits every element of $H$. Setting $e_\alpha(G)=Y$ defines the desired embedding
$e_\alpha:\B_\alpha\to\PNfin$

Finally, $e=\bigcup\{e_\alpha : \alpha <\omega_1\}$ is then the required embedding.
\end{proof}

\begin{theorem}[$\CH$]\label{C_aut_emb} Let $\B$ be a Boolean algebra of size $\c$  and $\mathfrak u(\B)>\omega$. 
Let $S\subseteq \Aut(\B) \setminus\{\id\}$, $|S|\le \c$ be given. 
There is an embedding $e: \B \to  \PNfin$ such that no element of $S$  can be extended to an element of  
$\Aut(\PNfin) $. 
\end{theorem}

\begin{proof} 
Let $\{h_\alpha  : \alpha <\omega_1\}$ be an enumeration of $S$ and let 
$H_\alpha $ be the homeomorphism of $\St(\B)$ corresponding to $h_\alpha $.    
Recursively choose ultrafilters $p_\alpha \in \St(\B)$  so that 
$\{p_\alpha : \alpha <\omega_1\} \cap \{H_\alpha (p_\alpha) : \alpha <\omega_1\}=\emptyset $. 
Then let $e: \B\to \PNfin$ be the embedding described in
Lemma \ref{split_emb} with respect to the set $\{p_\alpha :\alpha<\omega_1\}$.

Suppose that, for some $\alpha  <\omega_1$, 
$h_\alpha \in S$ has an extension $\tilde h_\alpha  \in
\Aut(\PNfin)$. As $C_\alpha <e(B)$ for every $B\in
p_\alpha $, $\tilde h_\alpha (C_\alpha )<\tilde
h_\alpha (e(B))=e(h_\alpha (B))$ for every $B\in p_\alpha
$. As $H_\alpha (p_\alpha )\ne p_\beta $ and $C_\alpha \meet C_\beta =\0$ for every
$\beta<\omega_1$, $\tilde h_\alpha (C_\alpha )\meet C_\beta =\0$ for all $\beta \in \omega_1$. 
Choose a $B\in \B$ which splits $\tilde h_\alpha (C_\alpha )$ and set $D=h_\alpha^{-1} (B)$. 
Note that $C_\alpha <e(D)$ or $C_\alpha <e(D)^c$ depending on whether $D\in p_\alpha$ or $D^c\in p_\alpha $.
In both cases this contradicts the fact that $e(B)$  splits $\tilde h_\alpha (C_\alpha )$.
\end{proof}

In particular, assuming $\CH$, every Boolean algebra $\B$ of size $\c$ with no ultrafilters of countable 
character having only $\c$ many automorphisms can be embedded into $\PNfin$ so that no non-trivial 
automorphism of $\B$ extends. Examples of such algebras include the measure algebra
or the algebra $\RO(2^\omega)$.\par
Aiming for a similar result for algebras with many automorphisms we first need the following lemma of independent interest.
Recall that the \emph{reaping number $\mathfrak r(\B)$} of a Boolean algebra $\B$ denotes the minimal size of a family $\calR\subseteq
\B^+$ such that for every $A\in\B^+$ there is an $R\in\calR$ such that $R\le B$ or $R\meet B=\0$ (in other words, no
$A\in\B^+$ splits all elements of $\calR$). A subset $D$ of a topological space $X$ is \emph{$G_\delta$-dense} if $D\cap G\neq\emptyset$
for every $G_\delta$ set $G\subseteq X$.

\begin{lemma}[$\CH$]\label{dense} Let $\B$ be a Boolean algebra of cardinality $\c$ and with $\mathfrak r (\B)>\omega$.
Then there is a $G_\delta$-dense $D\subseteq\St(\B)$ of size $\aleph_1$ such that every countable subset of $D$ 
is relatively discrete.
\end{lemma}

\begin{proof} Enumerate $\B$ as $\{B_\alpha:\alpha<\c\}$ and list all decreasing sequences in $\B^+$ as 
$\{\<B_{\alpha,n}\>_{n\in\omega}:\alpha<\c\}$ in such a way, so that $\{B_{\alpha,n}:n\in\omega\}\subseteq\{B_\beta:\beta<\alpha\}$ 
for every $\alpha<\c$. We will recursively construct countable families $\{\calF_{\alpha,\beta}:\beta<\alpha\}$ of countably generated 
filters on $\B$ so that for every $\alpha<\omega_1$

\begin{enumerate}
\item $(\forall \beta\leq\alpha)( B_\alpha\in \calF_{\alpha+1,\beta}\text{ or }  {B_\alpha}^c\in \calF_{\alpha+1,\beta})$,
\item $(\forall \gamma<\beta\leq\alpha)(\calF_{\beta,\gamma}\subseteq\calF_{\alpha,\gamma})$,
\item $(\exists\beta\leq\alpha)(\{B_{\alpha,n}:n\in\omega\}\subseteq \calF_{\alpha+1,\beta})$,
\item $(\forall \beta<\alpha)(\exists F_\beta\in\calF_{\alpha+1,\beta})(\beta\neq\gamma\Rightarrow 
F_\beta\meet F_\gamma=\0)$.
\end{enumerate} 

Assume that $B_0=\1$ and $B_1=\0$ and set $\calF_{1,0}=\{\1\}$. For limit $\alpha<\omega_1$ let 
$\calF_{\alpha,\beta}=\bigcup\{\calF_{\gamma,\beta}:\beta<\gamma<\alpha\}$. Given $\{\calF_{\alpha,\beta}:\beta<\alpha\}$, fix
a bijection $\varphi:\alpha\to\omega$ and construct
$\{\calF_{\alpha+1,\beta}:\beta<\alpha\}$ as follows:\par

Recursively pick $R_n$ ($n\in\omega$) so that
$R_n$ splits all elements of the family of all non-zero finite intersections from the family 
$\bigcup\{\calF_{\alpha,\beta}:\beta<\alpha\}\cup\{R_i:i<n\}$ and let 
$$F^\alpha_\beta= (\bigwedge_{m<\varphi(\beta)} R_m)\meet {R_{\varphi(\beta)}}^c.$$
Let, for all $\beta<\alpha$, $\calF_{\alpha+1,\beta}$ be the filter generated by $\calF_{\alpha,\beta}\cup\{F^\alpha_\beta, B_\alpha\}$ 
if possible, otherwise let  $\calF_{\alpha+1,\beta}$ be the filter generated by $\calF_{\alpha,\beta}\cup\{F^\alpha_\beta, {B_\alpha}^c\}$
Finally, choose $\calF_{\alpha+1,\alpha}$ so that 
$\{B_{\alpha,n}:n\in\omega\}\subseteq \calF_{\alpha+1,\alpha}$ and so that
$B_\beta\in \calF_{\alpha+1,\alpha}\text{ or }  {B_\beta}^c\in \calF_{\alpha+1,\beta}$ for all $\beta\leq\alpha$.
This concludes the recursive construction.

In the end set $x_\alpha=\bigcup_{\alpha<\beta}\calF_{\beta,\alpha}$ and note that:
\begin{itemize}
\item $x_\alpha$ is an ultrafilter,
\item the family $\{F^\alpha_\beta: \beta<\alpha\}$ witnesses that $\{x_\beta:\beta<\alpha\}$ is relatively discrete, and
\item $\{B_{\alpha,n}:n\in\omega\}\subseteq x_\alpha$ so $D=\{x_\alpha:\alpha<\omega_1\}$
is $G_\delta-$dense in $\St(\B)$.
\end{itemize}
\end{proof}

\begin{lemma}[$\CH$]\label{P-set} There are $\c$ many non-homeomorphic separable rigid P-sets in $\N^*$.
\end{lemma}

\begin{proof} Follows immediately from results of Dow, Gubbi and Szyma\' nski \cite{DGS}  and Balcar, Frankiewicz and Mills \cite{BFM}.
In \cite{DGS} the authors prove that there are $2^\c$ many pairwise non-homeomorphic rigid separable extremally
disconnected compact spaces. In \cite{BFM} it is shown that, assuming $\CH$, any compact zero-dimensional F-space of weight at most $\c$
is homeomorphic to a nowhere dense P-set in $\N^*$ (see Theorem 1.4.4 of \cite{vM}). As every extremally disconnected space is an F-space
and compact separable spaces are of weight less or equal to $\c$ the lemma easily follows.
\end{proof}

The proof of the next theorem is an advanced bookkeeping argument, moreover obscured by the fact that 
it uses both the language of Boolean algebra and the dual language of topology. In order not to cloud
the argument any further by unnecessary notation we identify (denote by the same symbol) an element of a Boolean
algebra and the corresponding clopen subset of the Stone space of the algebra.

\begin{theorem}[$\CH$]\label{main_bad_CH} Let $\B$ be a Boolean algebra of size $\c$ such that $\mathfrak r(\B)>\omega$. 
There is an embedding $e: \B \to  \PNfin$ such that no element of $\Aut(\B)$ other than the identity 
can be extended to an automorphism of $\PNfin $. 
\end{theorem}

\begin{proof} Enumerate $\B$ as $\{B_\alpha:\alpha<\omega_1\}$ and let $D$ be a 
$G_\delta-$dense subset of $\St(\B)$ of size $\aleph_1$ such that every countable subset of $D$ is relatively discrete.
We will construct the embedding $e:\B\to\PNfin$ as an increasing union of a chain of embeddings $\{e_\alpha:\alpha<\omega_1\}$,
$e_\alpha:\B_\alpha\to\PNfin$ where $\B_\alpha$ is a countable subalgebra of $\B$ containing $\{B_\beta:\beta<\alpha\}$. During the 
construction we will recursively choose points $x_\beta\in D$ together with an $A_\beta\in\PNfin$, split into $A_\beta^r$ and $A_\beta^s$
with the intention  that 
\begin{enumerate}
\item $E^{-1}(x_\beta)\subseteq A_\beta$,
\item $E^{-1}(x_\beta)\cap A_\beta^r$ is regular closed subset of $\omega^*$,
\item $E^{-1}(x_\beta)\cap A_\beta^s$ is a separable P-set,
\item $E^{-1}(x_\beta)\cap A_\beta^s$ and $E^{-1}(x_\gamma)\cap A_\gamma^s$ are not homeomorphic as long as $\beta\neq\gamma$,
\item $\bigcup_{\beta<\omega_1}E^{-1}(x_\beta)$ is dense in $\N^*$.
\end{enumerate}
where $E$ denotes the continuous surjection $\omega^*\to\St(\B)$ dual to $e$.\par

Assume first that this can be accomplished and let $h$ be an automorphism of the Boolean algebra $\B$ that can be extended
to an automorphism $\bar h$. 
Let $H$ denote the  autohomeomorphism of $St(\B)$ dual to $h$ and similarly,
$\bar H$ denotes the  autohomeomorphism of $\N^*$ dual to $\bar h$.
Note first that, for every $\alpha<\omega_1$ there is a $\beta<\omega_1$ 
such that $H(x_\alpha)=x_\beta$ as by (2) and (5) the points $x_\alpha$ are exactly those points $x\in \St(\B)$
for which  $E^{-1}(x)$ has non-empty interior. Then, however, $\bar H\restr E^{-1}(x_\alpha)$ is a homeomorphism
between $E^{-1}(x_\alpha)$ and $E^{-1}(x_\beta)$ and furthermore $\bar H[E^{-1}(x_\alpha)\cap A_\alpha^s] =
E^{-1}(x_\beta))\cap A_\beta^s$. By (4) this implies that $H(x_\alpha)=x_\alpha$ for all $\alpha<\omega_1$ and as 
$\{x_\alpha:\alpha<\omega_1\}$ is dense in $\St(\B)$ (by (5)), it follows that $h=\id$.
 
\bigskip

Next we will specify the promises which will ensure that the conditions 
(1)--(5) are satisfied.
To that end enumerate $\PNfin$ as $\{X_\alpha:\alpha<\omega_1\}$.
At stage $\alpha$ (when we build $\B_{\alpha+1}$ and $e_{\alpha+1}$,
 and choose $x_\alpha$) we look, for every $\beta<\alpha$, at 
the family $\{e_\alpha(B)\cap A^r_\beta\cap X_\alpha: B\in x_\beta\cap \B_\alpha\}$  and 
diagonalize it (if possible) by $R_{\alpha,\beta}\in\PNfin^+$. Otherwise, if $e_\alpha(B)\cap A^r_\beta\cap X_\alpha=\0$ 
for some $B\in x_\beta\cap \B_\alpha$, we will set $R_{\alpha,\beta}=\0$. Similarly, $R_{\beta ,\alpha}$ will be a diagonalization 
(if possible) of the family $\{e_\alpha(B)\cap A^r_\alpha\cap X_\beta: B\in x_\alpha\cap \B_\alpha\}$ for $\beta\leq\alpha$. 
The promise 
is that
$$
(\forall\beta\leq \alpha)(\forall B\in x_\beta)(R_{\alpha,\beta}< e(B))
\leqno(A)
$$
This will ensure that $E^{-1}(x_\beta)\cap A_\beta^r=\cl{\bigcup_{\beta\leq \alpha<\omega_1} R_{\alpha,\beta}}=R_\beta$ 
is regular closed.
We will use Lemma \ref{P-set} to choose a closed separable P-set $S_\alpha \subseteq A^s_\alpha$ (non-homeomorphic to any of 
the previous choices
$S_\beta$, $\beta<\alpha$, and disjoint from them) and its decreasing neighborhood base 
(a P-filter base on $\PNfin$) $\{T_{\alpha,\gamma}:\gamma\geq\alpha\}$ with $T_{\alpha,\alpha}=A^s_\alpha$ and promise that
$$(\forall \gamma\geq\alpha)(\forall B\in \B_\gamma\cap > x_\alpha)( S_\alpha\subseteq e(B)\cap A_\alpha^s)\text{ and }
\leqno(B)$$
$$
(\forall \gamma\geq\alpha)(\exists B\in\B_\gamma\cap x_\alpha)
(e(B)\cap A_\alpha^s\subseteq T_{\alpha,\gamma})
$$
to ensure that $E^{-1}(x_\alpha)\cap A_\alpha^s=S_\alpha$. Further, we promise
$$
(\exists \beta\leq\alpha)(S_\beta\cap X_\alpha\neq^*\emptyset
\text{ or }R_{\beta,\alpha}>\0)
\leqno(C)
$$
to guarantee that $\bigcup_{\beta<\omega_1}E^{-1}(x_\beta)$ is dense in $\omega^*$. Finally, in order to make sure that 
the recursive construction will not prematurely terminate we will require that
$$
\rng (e)\cap \< R_{\gamma,\beta}:\gamma,\beta<\alpha \>=\{\0,\1\},
\leqno(D)
$$
where $\< R_{\gamma,\beta}:\gamma,\beta<\alpha \>$ denotes the subalgebra of $\PNfin$ generated by 
$\{R_{\gamma,\beta}:\gamma,\beta<\alpha\}$.

It is easy to see that if we succeed in constructing algebras $\B_\alpha$ such that $\B=\bigcup_{\alpha<\omega_1}\B_\alpha$ and embeddings
$e_\alpha:\B_\alpha\to\PNfin$, together with choosing $x_\alpha\in D$, $A^s_\alpha$ and $A^r_\alpha$ as well as
$R_{\alpha,\beta}$ and $T_{\alpha,\beta}$ so that the above conditions are satisfied, then so are the conditions (1)--(5).

\bigskip

Assume that an increasing chain $\{\B_\beta :\beta\in \alpha \}$ of countable subalgebras of $\B$ and compatible
embeddings $e_\beta :\B_\beta \to \PNfin$ have already been constructed and set
$\C=\bigcup\{\B_\beta :\beta\in \alpha \}$ as well as $c=\bigcup\{e_\beta :\beta\in \alpha \}$. Let
$\C^+$ be the subalgebra of $B$ generated by  $\C\cup\{B_\alpha\}$. Assume that $B_\alpha\not\in \C$ and let
$$F=\{c(C):C<B_\alpha\}\cup\{ R_{\beta,\gamma}:\beta, \gamma<\alpha \ \& \ B_\alpha\in x_\beta\}\cup\{T_{\beta,\gamma_\beta}:
\beta<\alpha \ \& \ B_\alpha\in x_\beta\},$$
$$G=\{c(C):C>B_\alpha\}\cup\{ {R_{\beta,\gamma}}^c:\beta, \gamma<\alpha \ \& \  B_\alpha\not\in x_\beta\}\cup\{{T_{\beta,\gamma_\beta}}^c:
\beta<\alpha \ \&\  B_\alpha\not\in x_\beta\},$$
$$H=\{c(C):C\not<B_\alpha\ \& \  C\not>B_\alpha\},$$
where $\gamma_\beta$ is such that:
\begin{itemize}
\item $ (\forall B\in \C)(S_\beta\subseteq c(B)\rightarrow T_{\beta,\gamma_\beta} < c(B))$,
\item The family $\{ T_{\beta,\gamma_\beta}:\beta<\alpha\}$ is pairwise disjoint,
\item $(\forall B\in \C)(c(B)\not\in\<\{ R_{\gamma,\beta}:\beta,\gamma<\alpha\}\cup\{ T_{\beta,\gamma_\beta}:\beta<\alpha\}\>$.
\end{itemize}
It is now easy to check that $\calR_\omega$ applies (here is where
the technical requirement (D) comes in handy), and $c$ can be extended to an embedding $c^+:\C^+\to\PNfin$.

For every $\beta<\alpha$ consider the family $\{c^+(B)\cap A^r_\beta\cap X_\alpha: B\in x_\beta\cap \C^+\}$.  If possible,
diagonalize it by $R_{\alpha,\beta}\in\PNfin^+$. Otherwise, fix $C_\beta\in x_\beta\cap \C^+$ such that 
$c^+(C_\beta)\cap A^r_\beta\cap X_\alpha=\0$ and set $R_{\alpha,\beta}=\0$.
If there is a $\beta <\alpha$ such that $R_{\alpha,\beta}>\0$ or such that $X_\alpha\cap S_\beta\neq\emptyset$ let 
$\{D_\beta:\beta<\alpha\}$ be a pairwise disjoint family in $\B$ such that $D_\beta\in x_\beta$ (recall that all
countable subsets of $D$ are relatively discrete). If $R_{\alpha,\beta}=\0$ and $X_\alpha\cap S_\beta=\emptyset$ for all  $\beta <\alpha$, 
require moreover that $D_\beta\leq C_\beta$. By picking the $D_\beta$'s sufficiently small
one can make sure that there is a $D_\alpha\in\B^+$
disjoint from all $D_\beta$, $\beta<\alpha$. Let $\C^{++}$ be the subalgebra of $\B$
generated by $\C^+\cup\{D_\beta:\beta<\alpha\}$. Enumerate $\alpha$ as $\{\beta_n:n\in\omega\}$ and extend $c^+$ to $c^{++}:\C^{++}\to\PNfin$
step by step, each time determining $c^{++}(D_{\beta_n})$ so that
\begin{itemize}
\item $(\forall \gamma\leq\alpha)(R_{\beta_n,\gamma}<c^{++}(D_{\beta_n}))$,
\item $T_{\beta_n,\gamma_{\beta_n}}<c^{++}(D_{\beta_n})$,
\item $(\forall i< n)(c^{++}(D_{\beta_i})\meet c^{++}(D_{\beta_n})=\0)$,
\item $(\forall i\neq n)(\forall\gamma\leq\alpha) (R_{\beta_i,\gamma}\meet c^{++}(D_{\beta_n})=\0)$,
\item $(\forall i\neq n)(\forall\gamma\leq\alpha) (T_{\beta_i,\gamma_{\beta_i}}\meet c^{++}(D_{\beta_n})=\0)$,
\item $(\forall B\in\<\C^+\cup\{D_{\beta_i}:i<n\}\>)(D_{\beta_n}<B \ \rightarrow c^{++}(D_{\beta_n})<c^{++}(B))$,
\item $c^{++}(D_{\beta_n})\meet A^s_{\beta_n}\leq T_{\beta_n,\alpha}$,
\item $\rng(c^{++})\cap \<R_{\beta,\gamma}:\beta,\gamma\leq\alpha\>=
       \{\0,\1\}$.
\end{itemize}
These are all compatible demands. Pick an ultrafilter $u$ on $\C^{++}$ containing ${D_\beta}^c$ for all $\beta<\alpha$ such that
the family $\{X_\alpha\}\cup\{c^{++}(B):B\in u\}$ is centered. Choose $A_\alpha$ below $X_\alpha\meet c^{++}(B)$ for $B\in u$ 
and split it into two pieces $A^s_\alpha$ and $A^r_\alpha$. Use Lemma \ref{P-set} to choose a closed separable P-set $S_\alpha \subseteq A^s_\alpha$
non-homeomorphic to any of the previous choices $S_\beta$, $\beta<\alpha$ and let $\{T_{\alpha,\gamma}:\gamma\geq\alpha\}$ with 
$T_{\alpha,\alpha}=A^s_\alpha$ be a decreasing neighborhood base of $S_\alpha$. Finally, pick  $x_\alpha\in D$ extending $u$ (this can be done
as $D$ is $G_\delta$-dense) and determine $R_{\alpha,\beta}$ for $\beta\leq\alpha$.\par
This finishes the construction. It is not difficult to verify that the promises
 (A)--(D) are fulfilled, the $D_\beta$'s being 
the witnesses required by (B).
\end{proof}

\section{Concluding remarks}

In this section we present some examples and open questions.

\begin{defn} Given an embedding $e:\B\to\A$ of Boolean algebras, let 
\item{} $E_e=\{h\in\Aut(\B): (\exists \bar h\in \Aut(\A))( e\circ h=\bar h\circ e)\}$,
\item{} $R_e=\{h\in\Aut(\A): h\restr \rng(e)\in\Aut(\rng(e))\}$,
\item{} $K_e=\{h\in\Aut(\A): h\restr \rng(e)= \id\}$.
\end{defn}

All three of these sets are easily seen to be groups, moreover,  $K_e$ is a normal subgroup of
$R_e$ and $E_e$ is isomorphic to the quotient group
$R_e/K_e$. In this notation, if $e:\B\to \PNfin$ is the
embedding described in Theorem \ref{grzech}, then $K_e$ is a direct summand
of $E_e$ and therefore  $R_e$ is isomorphic to $E_e\oplus K_e$. A
similar decomposition is trivially realized when $e$ is a dense
embedding, as in this case $K_e=\{\id\}$.\par
However, even in the particular case of $\PNfin$,  it is not immediately obvious whether any
non-dense embedding leads to a trivial $K_e$. 
Clearly, the range of any such $e$ has
to form a splitting family in $\PNfin$. We will show next that there is (in ZFC) an embedding $e$
of the free Boolean algebra on $\c$ generators into $\PNfin$ such that 
the only automorphism of $\PNfin$ which restricts to an automorphism of $\CO(2^\c)$ is the identity.\par
Recall that a family $\calI$ of subsets of a set $X$ is \emph{independent} if 
$\bigcap \calF_0\setminus \bigcup \calF_1\neq\emptyset$ for all pairs $\calF_0, \calF_1$
of disjoint finite subsets of $\calI$. It is well-known that there is an independent family of size $\c$.

\begin{theorem}\label{free_emb} There is an embedding $e: \CO(2^\c)\to \PNfin$ such that $K_e=\{\id\}$.
\end{theorem}

\begin{proof} Let $\{I_\alpha:\alpha<\c\}$ be an independent family of subsets of $\N$. Enumerate
all pairs of disjoint nowhere dense subsets of the rational numbers $\Q$ as $\{\<M_\alpha,N_\alpha\>:\alpha<\c\}$.
Set
$$J_\alpha=\left(M_\alpha\cup\bigcup_{n\in I_\alpha} [n,n+1)\cap\Q\right)\setminus N_\alpha$$
for every $\alpha<\c$. Note that $\{J_\alpha:\alpha<\c\}$ is an independent family of subsets of $\Q$.
Identify $\PNfin$ with $\calP( \Q) /\text{fin}$ and set
$$e(\{f\in 2^\c: f(\alpha)=1\})= [J_\alpha],$$
defining thus an embedding $e: \CO(2^\c)\to \calP( \Q) /\text{fin}$.
All that remains to be checked is that $K_e=\{\id\}$. To that end let $h\in\Aut(\calP( \Q) /\text{fin})\setminus\{\id\}$, i.e., there is
an infinite $A\subseteq \Q$ such that $h([A])\neq [A]$. It is easy to find an infinite nowhere dense subset of $A$, say $M$, such that there 
is an
$N\in h([M])$ nowhere dense, disjoint from $M$. The pair $\<M,N\>$ was enumerated as $\<M_\alpha,N_\alpha\>$ for some $\alpha<\c$. 
Then, however,
$$h([J_\alpha])\geq h([M_\alpha])=[N_\alpha]\text{ and } [J_\alpha]\meet [N_\alpha]=\0$$
so $h([J_\alpha])\neq [J_\alpha]$ and $h\not\in K_e$.
\end{proof}

One might wonder to what extent are the very strong assumptions in Theorem \ref{goodemb_MA}   
necessary. For instance, Theorem \ref{goodemb_MA} only applies to Boolean algebras satisfying the countable chain condition
for if a Boolean algebra $\B$ admits an embedding which lifts, it has to be c.c.c. 
The following example seems to indicate that one can not hope for a substantial weakening
of the assumptions. Let $\A(\kappa)$ denote the Boolean algebra of finite and co-finite subsets of $\kappa$ or
equivalently the algebra of clopen subsets of the one point compactification of a discrete space of size $\kappa$,
i.e., the simplest Boolean algebra having an antichain of size $\kappa$.\par
Recall that a family $\calA\subseteq \calP(\N)$ is \emph{almost disjoint} if every two distinct members of $\calA$
have finite intersection. An almost disjoint family $\calA$ of size $\aleph_1$  
is \emph{Luzin} if there is no $B\subseteq\N$ such that $|\{A\in\calA: |A\setminus B|\text{ is finite}\}|=
|\{A\in\calA: A\cap B\text{ is finite}\}|=\aleph_1$. Inspired by the ingenious construction of Hausdorff, N. Luzin
(\cite{Lu}) showed that there are Luzin almost disjoint families in ZFC.

\begin{proposition} There is an embedding $e:\A(\omega_1)\to\PNfin$ such that $E_e\neq \Aut(\A(\omega_1))$.
\end{proposition}

\begin{proof} Let $\calA$ be a Luzin almost disjoint family of size $\aleph_1$. Identify $\N$ with $\N\times 3$
and let, for $X\subseteq\N$ and $i\in 3$, $X_i=X\times\{i\}$. Let $\A$ be the subalgebra of $\P(\N\times 3)/\fin$
generated by $[A_i]$, for $A\in\calA$ and $i\in 3$. $\A$ is isomorphic to $\A(\omega_1)$.
Let $\varphi$ be a bijection between $\{A_1:A\in\calA\}\cup\{A_2:A\in\calA\}$ and $\{A_0:A\in\calA\}$
and let $h$ be the induced automorphism of $\calA$. Assume that $h$ has an extension $H\in\Aut(\P(\N\times 3)/\fin)$ and consider
$H([\N_1])$. The sets $\{A\in\calA: [A_0]\leq H([\N_1])\}$
and $\{A\in\calA: [A_0]\meet H([\N_1])=\0\}=\{A\in\calA: [A_0]\leq H([\N_2])\}$ are both uncountable
contradicting that $\calA$ was a Luzin family.
\end{proof}

A natural question is whether the assumption in Theorem \ref{main_bad_CH} can be weakened to $\mathfrak u(\B)>\omega$, which
would, of course make, Theorem \ref{C_aut_emb} obsolete. 

\begin{question} Assume $\CH$. Let $|\B|=\mathfrak u(\B)=\omega_1$. Is there an embedding $e: \B \to  \PNfin$ such that no 
element of $\Aut(\B)$ other than the identity can be extended to an automorphism of $\PNfin $?
\end{question}

Note that for the embeddings constructed in Theorem \ref{C_aut_emb} and Theorem \ref{main_bad_CH}, $E_e=\{\id\}$ but
$R_e\not =\{\id\}$. This suggests the following question:

\begin{question} Assume $\CH$. For which Boolean algebras $\B$ is there an embedding $e: \B \to  \PNfin$ such that no 
element of $\Aut(\PNfin)$ other than the identity restricts to an automorphism of $\B$?
\end{question}

To this question we have a very partial answer.

\begin{theorem}[$\CH$]\label{R_e} There is an embedding $e:\PNfin\to\PNfin$ such that $R_e=\{\id\}$.
\end{theorem}

\begin{proof} Enumerate $\PNfin$ as $\{B_\alpha:\alpha<\omega_1\}$ and let $p\in\N^*$ be a P-point. Let
$T(p)$ be the type of $p$, i.e., the set of those $q\in\N^*$ which can be sent to $p$ by a permutation of $\N$.
Enumerate $T(p)$ as $\{p_\alpha:\alpha<\omega_1\}$ and note that it is a  
$G_\delta-$dense subset of $\N^*$ every countable subset of which is relatively discrete.
List as $\{X_\alpha, Y_\alpha, \calX_\alpha, \calY_\alpha\}$ all quadruples such that
\begin{itemize}
\item $X_\alpha, Y_\alpha\in \PNfin$, $X_\alpha\subseteq Y_\alpha$,
\item $\calX_\alpha, \calY_\alpha\in [\PNfin]^\omega$,  $\calX_\alpha\cup\calY_\alpha$ is a pairwise disjoint family,
\item $(\forall Y\in\calY_\alpha)(Y\cap X_\alpha=\0)$,
\item $(\forall X\in\calX_\alpha)(X\subseteq Y_\alpha)$.
\end{itemize}
As before we  will construct the embedding $e:\PNfin\to\PNfin$ as an increasing union of a chain of embeddings $\{e_\alpha:\alpha<\omega_1\}$,
$e_\alpha:\B_\alpha\to\PNfin$ where $\B_\alpha$ is a countable subalgebra of $\B$ containing $\{B_\beta:\beta<\alpha\}$.\par 
At stage $\alpha<\omega_1$ of the construction we will use Lemma \ref{P-set} to
recursively choose a rigid separable P-set $S_\alpha \subseteq \N^*$ 
non-homeomorphic to and disjoint from all of the previous choices $S_\beta$, $\beta<\alpha$ and its decreasing neighborhood 
base $\{T_{\alpha,\gamma}:\gamma\geq\alpha\}$, and $D_\alpha\in \PNfin^+$ 
promising that:
\begin{enumerate}
\item $(\forall \gamma\geq\alpha)(\forall B\in \B_\gamma\cap p_\alpha)( S_\alpha\subseteq e_\gamma(B))$,
\item $(\forall \gamma\geq\alpha)(\exists B\in \B_\gamma\cap p_\alpha)( e_\gamma(B)\subseteq T_{\alpha,\gamma})$,
\item $(\forall \beta\leq\alpha)(\forall \gamma\ge\alpha)(B_\beta$ does not split $S_\gamma)$,
\item $X_\alpha\cup\bigcup \{ X:X\in\calX\}\subseteq e_\alpha(D_\alpha)\subseteq Y_\alpha\setminus\bigcup\{ Y:Y\in\calY\}$,
\item $(\exists\gamma\geq\alpha)( S_\alpha\cap B_\alpha\neq\0)$.
\end{enumerate} 
This can be accomplished in a very similar way  to the appropriate part of the proof of Theorem \ref{main_bad_CH}.
Note the double role the $B_\alpha$'s play.
Assume that the construction has been successfully completed and let $E:\N^*\to\N^*$ be the dual map to the embedding 
$e=\bigcup_{\alpha<\omega_1}e_\alpha$.
\begin{itemize}
\item $(\forall\alpha<\omega_1)(E^{-1}(p_\alpha)=S_\alpha)$
\end{itemize}
This follows easily from (1) and (2). Now assume that $q\in\N^*$ is a P-point of a different type that
$p$. We will show that
\begin{itemize}
\item $|E^{-1}(q)|=1$
\end{itemize}
Aiming for a contradiction assume that $q=E(r)=E(s)$ for some $r\neq s\in\N^*$. Let $Y\in s$ be a neighborhood
of $s$ disjoint from r. By (3) the set $\calA=\{\beta: S_\beta$ is split by $Y\}$ is countable. As $q$ is a P-point 
there is a $Q\in q$ disjoint from $\{p_\beta:\beta\in \calA\}$. Then $X=Y\cap e(Q)$ is a
neighborhood of $s$ such that $X\cap S_\beta=\emptyset$ for every $\beta\in\calA$. Let $\calB= 
\{\beta: S_\beta$ is split by $X\}$. Clearly, $\calB$ is also countable. As $\calA\cup\calB$ is a pairwise disjoint collection
of P-sets, there is a ``swelling''  of $\calA$ and $\calB$ respectively to families $\calY$ and $\calX$ of pairwise disjoint clopen sets
such that  $Y'\cap X=\0$ for every $Y'\in\calY$ and $X'\subseteq Y$ for every $X'\in\calX$. Now, the quadruple
$\{X, Y, \calX, \calY\}$ was enumerated as $\{X_\alpha, Y_\alpha, \calX_\alpha, \calY_\alpha\}$ for some $\alpha<\omega_1$.
Consider $D_\alpha$. By (4) $e(D_\alpha)\in s\setminus r$ and so $E(s)\neq E(r)$.\par
Let $h\in R_e$ and let $H$ be the dual autohomeomorphism of $\N^*$, let $\tilde h$ be the restriction of $h$ onto the
range of $e$ and let $\tilde H$ be the map dual to $\tilde h$. Note first that 
any homeomorphism sends P-points to P-points. In particular, $\tilde H(p_\alpha)$ has to be a P-point, and moreover,
$H\restr E^{-1}(p_\alpha)$ has to be a homeomorphism between   $ E^{-1}(p_\alpha)$ and $E^{-1}(\tilde H(p_\alpha))$. It follows
that $\tilde H(p_\alpha)= p_\alpha$ for every $\alpha<\omega_1$ as $ E^{-1}(p)$ is not homeomorphic to $ E^{-1}(p_\alpha)$
for any other P-point in $\N^*$. By density of $\{p_\alpha:\alpha<\omega_1\}=T(p)$, $\tilde H=\id$ and consequently, 
$H(r)=s$ implies that $E(r)=E(s)$. To conclude the argument recall that all sets $S_\alpha=E^{-1}(p_\alpha)$ were rigid,
hence $H\restr E^{-1}(p_\alpha)=\id$ for every $\alpha<\omega_1$. By (6), $\bigcup_{\alpha<\omega_1}E^{-1}(p_\alpha)$ is dense
in $\N^*$ so $H=\id$.
\end{proof}

As with most results contained here the situation becomes quite different if the assumption of the Continuum Hypothesis
is dropped. One might have problems even formulating reasonable conjectures for at least two distinct reasons:
(1) There may not be an  embedding of the Boolean algebra into $\PNfin$ \cite{DH}, and (2) there may not be enough
automorphisms of $\PNfin$ \cite{SS}.
For instance, assuming the Proper Forcing Axiom, Theorem \ref{R_e} no longer holds:

\begin{theorem}[PFA] 
$E_e\not =\{\id\}$ for every embedding $e:\PNfin\to\PNfin$.
\end{theorem}

\begin{proof} 
Follows immediately from a result of I. Farah (see \cite{Fa}, Theorem 3.8.1). He has shown that, assuming $\PFA$
 if $e:\PNfin\to\PNfin$ is an embedding
then there is an $A\in [\N]^\omega$ and a one-to-one map $f:A\to\N$ such that $e([B])=[f[B]]$ for every $B\subseteq A$.
\end{proof}

%\bibliography{names,longnames,research}
%\bibliographystyle{amsplain}

\providecommand{\bysame}{\leavevmode\hbox to3em{\hrulefill}\thinspace}
\providecommand{\MR}{\relax\ifhmode\unskip\space\fi MR }
\providecommand{\MRhref}[2]{%
  \href{http://www.ams.org/mathscinet-getitem?mr=#1}{#2}
}
\providecommand{\href}[2]{#2}

\end{document}